\documentclass[a4paper,oneside,10pt]{article} 

\title{ARTICLE - INEQUALITY ENTROPY}
\author{Rapha\"el Cerf & Matthias Gorny}
\date{\textit{MAJ: October 17th, 2014}}

\usepackage[english]{babel}
\usepackage{amsfonts}
\usepackage{cite}
\usepackage{amsthm}
\usepackage{amssymb}
\usepackage{amsmath}
\usepackage{dsfont}

\newcommand{\R}{\mathbb{R}}

\newcommand{\Ll}{\mathrm{L}}
\newcommand{\Fc}{\mathcal{F}}
\renewcommand{\r}{\rho}
\renewcommand{\l}{\lambda}

\newtheorem*{theo}{Theorem}
\newtheorem*{prop}{Proposition}

\begin{document}

\renewcommand{\contentsname}{Contents}
\renewcommand{\refname}{\textbf{References}}
\renewcommand{\abstractname}{Abstract}

\begin{center}
\begin{Huge}
A Lower Bound on the\medskip

Relative Entropy with Respect to\bigskip

a Symmetric Probability
\end{Huge}\bigskip \bigskip \bigskip \bigskip

\begin{Large} Rapha\"el Cerf \end{Large} \smallskip

\begin{large} {\it ENS Paris} \end{large}\bigskip \medskip

and \bigskip\medskip

\begin{Large} Matthias Gorny \end{Large} \smallskip
 
\begin{large} {\it Universit\'e Paris Sud and ENS Paris} \end{large} \bigskip \bigskip
\end{center}
\bigskip \bigskip \bigskip

\begin{abstract}
\noindent Let~$\rho$ and~$\mu$ be two probability measures on~$\mathbb{R}$ which are not the Dirac mass at~$0$. We denote by~$H(\mu|\rho)$ the relative entropy of~$\mu$ with respect to~$\rho$. We prove that, if~$\rho$ is symmetric and~$\mu$ has a finite first moment, then
\[ H(\mu|\rho)\geq \frac{\displaystyle{\left(\int_{\mathbb{R}}z\,d\mu(z)\right)^2}}{\displaystyle{2\int_{\mathbb{R}}z^2\,d\mu(z)}}\,,\]
with equality if and only if~$\mu=\rho$.
\end{abstract}
\bigskip \bigskip \bigskip \bigskip \bigskip



\bigskip\bigskip \bigskip\bigskip \bigskip

\newpage

\section{Introduction}

\noindent Given two probability measures~$\mu$ and~$\r$ on~$\R$, the relative entropy of~$\mu$ with respect to~$\r$ (or the Kullback-Leibler divergence of~$\r$ from~$\mu$) is
\[ H(\mu|\r)=\left\{\begin{array}{cl}
\displaystyle{\int_{\R} f(z)\ln f(z)\,d\r(z)} & \quad\mbox{if } \mu \ll \r  \mbox{ and } \displaystyle{f=\frac{d\mu}{d\r}}\\[0.3cm]
+\infty & \quad\mbox{otherwise}\,,
\end{array}\right.\]
where~$d\mu/d\r$ denotes the Radon-Nikodym derivative of~$\mu$ with respect to~$\r$ when it exists. In this paper, we prove the following theorem:

\begin{theo} Let~$\r$ and~$\mu$ be two probability measures on~$\R$ which are not the Dirac mass at~$0$. We suppose that
\[\int_{\R}|z|\,d\mu(z)<+\infty\,.\]
If~$\r$ is symmetric then
\[ H(\mu|\r)\geq \frac{\displaystyle{\left(\int_{\R}z\,d\mu(z)\right)^2}}{\displaystyle{2\int_{\R}z^2\,d\mu(z)}}\,,\]
with equality if and only if~$\mu=\r$.
\label{MainInequality}
\end{theo}

\noindent A remarkable feature of this inequality is that the lower bound does not depend on the symmetric probability measure~$\r$. We found the following related inequality in the literature (see lemma~3.10~of~\cite{concentrationBLM}): if $\r$ is a probability measure on~$\R$ whose first moment $m$ exists and such that
\[\exists v>0\qquad \forall \l\in \R\qquad\int_{\R}\exp(\l(z-m))\,d\r(z)\leq \exp\left(\frac{v\l^2}{2}\right)\,,\]
then, for any probability measure $\mu$ on~$\R$ having a first moment, we have
\[H(\mu|\r)\geq \frac{1}{2v}\left(\int_{\R}z\,d\mu(z)-m\right)^2\,.\]
Our inequality does not require an integrability condition. Instead we assume that $\r$ is symmetric.\medskip

\noindent The proof of the theorem is given in the following section. It consists in relating the relative entropy $H(\,\cdot\,|\r)$ and the Cram\'er transform~$I$ of $(Z,Z^2)$ when~$Z$ is a random variable with distribution~$\r$. We then use an inequality on~$I$ which we proved
initially in~\cite{CerfGorny}. We give here a simplified proof of
this inequality.

\section{Proof of the theorem}

\noindent Let~$\r$ and~$\mu$ be two probability measures on~$\R$ which are not the Dirac mass at~$0$. We first recall that $H(\mu|\r)\geq 0$, with equality if and only if $\mu=\r$.\medskip

\noindent We assume that $\r$ is symmetric and that $\mu$ has a finite first moment. We denote
\[\Fc(\mu)=\frac{\displaystyle{\left(\int_{\R}z\,d\mu(z)\right)^2}}{\displaystyle{2\int_{\R}z^2\,d\mu(z)}}\,.\]
If $\mu=\r$ then $\Fc(\mu)=0=H(\mu|\r)$. 
From now onwards we suppose that~$\mu\neq \r$. If the first moment of~$\mu$ vanishes or if its second moment is infinite, then~$\Fc(\mu)=0<H(\mu|\r)$. Finally, if~$\mu$ is such that~$H(\mu|\r)=+\infty$, then Jensen's inequality implies that
\[\Fc(\mu)\leq 1/2<H(\mu|\r).\]
In the following, we suppose that
\[\int_{\R}z \,d\mu(z)\neq 0,\qquad \int_{\R}z^2\,d\mu(z)<+\infty\,,\]
and that $H(\mu|\r)<+\infty$. This implies that~$\mu \ll \r$ and we set~$f=d\mu/d\r$. It follows from Jensen's inequality that, for any~$\mu$-integrable function~$\Phi$, 
\[\int_{\R}\Phi\,d\mu-H(\mu|\r)=\int_{\R}\ln\left(\frac{e^{\Phi}}{f}\right)\,d\mu\leq \ln \int_{\R}\frac{e^{\Phi}}{f}\,d\mu=\ln\int_{\R}e^{\Phi}\,d\r\,.\]
As a consequence
\[ \sup_{\Phi \in \Ll^1(\mu)}\,\left\{\,\int_{\R}\Phi\,d\mu-\ln\int_{\R}e^{\Phi}\,d\r\,\right\} \leq H(\mu|\r)\,.\]
In order to make appear the first and second moments of~$\r$, we consider functions~$\Phi$ of the form~$z\longmapsto uz+vz^2$,~$(u,v)\in \R^2$. 
This way we obtain
\[I\left(\int_{\R}z\,d\mu(z),\int_{\R}z^2\,d\mu(z)\right)\leq H(\mu|\r)\,,\]
where
\[\forall (x,y)\in \R^2 \qquad I(x,y)=\sup_{(u,v)\in \R^2}\,\left\{\,ux+vy-\ln\int_{\R}e^{uz+vz^2}\,d\r(z)\,\right\}\,.\]
The function $I$ is the Cram\'er transform of~$(Z,Z^2)$ when~$Z$ is a random variable with distribution~$\r$. In our paper dealing with a Curie-Weiss model of self-organized criticality~\cite{CerfGorny}, we 
proved with the help of the following inequality that, under some integrability condition,  
the function $(x,y)\longmapsto I(x,y)-x^2/(2y)$ has a unique global minimum on $\R\times\,]0,+\infty[$ at $\big(0,\int x^2\,d\r(x)\big)$.

\begin{prop} If $\r$ is a symmetric probability measure which is not the Dirac mass at~$0$, then
\[\forall x\neq 0\quad \forall y\neq 0\qquad I(x,y)> \frac{x^2}{2y}.\]
\end{prop}

\noindent We present here a proof of this proposition which is simpler than in~\cite{CerfGorny}.\medskip

\noindent{\bf Proof.} Let $x\neq 0$ and $y\neq 0$. By definition of $I(x,y)$, we have
\begin{align*}
I(x,y)&\geq x\times \frac{x}{y}+y\times \left(-\frac{x^2}{2y^2}\right)-\ln\int_{\R}\exp\left(\frac{xz}{y}-\frac{x^2z^2}{2y^2}\right)\,d\r(z)\\
&=\frac{x^2}{2y}-\ln\int_{\R}\exp\left(\frac{xz}{y}-\frac{x^2z^2}{2y^2}\right)\,d\r(z)\,.
\end{align*}
Let $(s,t)\in \R^2$. By using the symmetry of $\r$, we obtain
\begin{multline*}
\int_{\R} \exp(sz-tz^2)\,d\r(z)=\int_{\R}\exp(-sz-tz^2)\,d\r(z)\\
=\frac{1}{2}\left(\int_{\R} \exp(sz-tz^2)\,d\r(z)+\int_{\R}\exp(-sz-tz^2)\,d\r(z)\right)\\
=\int_{\R} \mathrm{cosh}(sz)\,\exp(-tz^2)\,d\r(z)\,.
\end{multline*}
We choose now $t=s^2/2$. We have the inequality
\[\forall u\in \R\backslash\{0\} \qquad
\mathrm{cosh}(u)\,\exp\left(-u^2/2\right)< 1\,.\]
Since $\r$ is not the Dirac mass at~$0$, the above inequality implies that
\[\forall s\neq 0\qquad\int_{\R} \mathrm{cosh}(sz)\,\exp\left(-\frac{s^2z^2}{2}\right)\,d\r(z)<1\,.\]
We finally choose $s=x/y$ and we get
\[\int_{\R}\exp\left(\frac{xz}{y}-\frac{x^2z^2}{2y^2}\right)\,d\r(z)<1\,.\]
As a consequence
\[I(x,y)\geq \frac{x^2}{2y}-\ln\int_{\R}\exp\left(\frac{xz}{y}-\frac{x^2z^2}{2y^2}\right)\,d\r(z)>\frac{x^2}{2y}\,,\]
which is the desired inequality.\qed\medskip

\noindent By applying the above proposition with
\[x=\int_{\R}z \,d\mu(z)\neq 0,\qquad y=\int_{\R}z^2\,d\mu(z)\in \,]0,+\infty[\,,\]
we obtain 
\[H(\mu|\r)\geq I\left(\int_{\R}z\,d\mu(z),\int_{\R}z^2\,d\mu(z)\right)>\Fc(\mu)\,.\]
This ends the proof of theorem~\ref{MainInequality}.

\bibliographystyle{plain}
\bibliography{biblio}
\addcontentsline{toc}{section}{References}
\markboth{\uppercase{References}}{\uppercase{References}}

\end{document}